\documentclass{article}

\usepackage{graphicx}
\usepackage{amsmath}
\usepackage{amssymb}
\usepackage{color}
\usepackage{hyperref}
\usepackage{epsfig}
\usepackage{tikz}
\usepackage{pgfplots}
\usepackage{caption}
\usepackage{float}
\allowdisplaybreaks

\newtheorem{theo}{Theorem}

\newtheorem{theorem}{Theorem}
\newtheorem{conj}{Conjecture}

\newtheorem{defi}{Definition}

\newtheorem{prop}{Proposition}

\newcommand{\NN}{{\mathbb N}}

\newcommand{\ZZ}{{\mathbb Z}}

\newcommand{\notiz}[1]{}

\title{Pair correlation of sequences $(\lbrace a_n \alpha \rbrace)_{n \in \NN}$ with maximal order of additive energy}
\author{Gerhard Larcher \footnote{The author is supported by the Austrian Science Fund (FWF), Project F5507-N26, which is a part of the Special Research Program “Quasi-Monte Carlo Methods: Theory and Applications” and Project I1751-N26.} \ , Wolfgang Stockinger \footnote{The author is supported by the Austrian Science Fund (FWF), Project F5507-N26, which is a part of the Special Research Program “Quasi-Monte Carlo Methods: Theory and Applications”.}}
\date{}
\begin{document}
\maketitle 
\begin{abstract}
We show for sequences $\left(a_{n}\right)_{n \in \NN}$ of distinct positive integers with maximal order of additive energy, that the sequence $\left(\left\{a_{n} \alpha\right\}\right)_{n \in \NN}$ does not have Poissonian pair correlations for any $\alpha$. This result essentially sharpens a result obtained by J.\ Bourgain on this topic.   
\end{abstract} 
\section{Introduction and statement of main results}

\begin{defi}
Let $\left\|\cdot\right\|$ denote the distance to the nearest integer. A sequence $\left(x_{n}\right)_{n \in \NN}$ in $\left[\left.0,1\right.\right)$ is said to have (asymptotically) Poissonian pair correlations, if for each $s > 0$ the pair correlation function
$$
R_{2} \left(\left[-s,s\right], \left(x_{n}\right)_{n}, N\right):= \frac{1}{N} \# \left\{1 \leq i \neq j \leq N \left| \left\|x_{i} -x_{j} \right\|\leq \frac{s}{N} \right.\right\} 
$$
tends to $2s$ as $N \rightarrow \infty$.
\end{defi}
It is known that if a sequence $\left(x_{n}\right)_{n \in \NN}$ has Poissonian pair correlations, then it is uniformly distributed modulo 1, cf., \cite{not4, not8, not15}. The converse is not true in general.\\

The study of Poissonian pair correlations of sequences, especially of sequences of the form $\left(\left\{a_{n} \alpha\right\}\right)_{n \in \NN}$, where $\alpha$ is irrational, and $\left(a_{n}\right)_{n \in \NN}$ is a sequence of distinct positive integers, is primarily motivated by certain questions in quantum physics, especially in connection with the Berry-Tabor conjecture in quantum mechanics, cf., \cite{not2, not11}. The investigation of Poissonian pair correlations was started by Rudnick, Sarnak and Zaharescu, cf., \cite{not12, not13, not14}, and was continued by many authors in the subsequent, cf., \cite{not3} and the references given there.\\

A quite general result which connects Poissonian pair correlations of sequences $\left(\left\{a_{n} \alpha\right\}\right)_{n \in \NN}$ to concepts from additive combinatorics was given in \cite{not3}:

For a finite set $A$ of reals the {\em additive energy} $E(A)$ is defined as
$$
E\left(A\right) := \sum_{a+b=c+d} 1,
$$
where the sum is extended over all quadruples $\left(a,b,c,d\right) \in A^{4}$. Trivially one has the estimate $\left|A^{2}\right| \leq E\left(A\right)\leq \left|A\right|^{3}$, assuming that the elements of $A$ are distinct. The additive energy of sequences has been extensively studied in the additive combinatorics literature, cf., \cite{not16}. In \cite{not3} the following was shown: 
\begin{theo}[in \cite{not3}]
Let $\left( a_n \right)_{n \in \NN}$ be a sequence of distinct integers, and let $A_{N}$ denote the first $N$ elements of this sequence. If there exists a fixed $\varepsilon > 0$ such that
$$
E \left(A_{N}\right) = \mathcal{O} \left(N^{3-\varepsilon}\right),
$$
then for almost all $\alpha$ the sequence $\left(\left\{a_{n} \alpha \right\}\right)_{n \in \NN}$ has Poissonian pair correlations.
\end{theo}
On the other hand Bourgain showed in \cite{not3} the following negative result:
\begin{theo}[in \cite{not3}]
If $E\left(A_{N}\right) = \Omega \left(N^{3}\right)$, then there exists a subset of $\left[0,1\right]$ of positive measure such that for \textbf{every} $\alpha$ from this set the pair correlations of $\left(\left\{a_{n} \alpha\right\}\right)_{n \in \NN}$ are \textbf{not} Poissonian.
\end{theo}
In \cite{not9}, the authors gave a sharper version of the result of Bourgain by showing that the set of exceptional values $\alpha$ from Theorem B has \textbf{full measure}.\\

It is the aim of this paper to show the best possible version of a result in this direction, namely: 
\begin{theorem}
If $E(A_N) = \Omega(N^3)$, then for \textbf{every} $\alpha$ the pair correlations of $(\lbrace a_n \alpha \rbrace)_{n \in \NN}$ are not Poissonian. 
\end{theorem} 
In fact, we conjecture that even more is true: 
\begin{conj}
If for almost all $\alpha$ the pair correlations of $(\lbrace a_n \alpha \rbrace)_{n \in \NN}$ are not Poissonian, then the pair correlations of this sequence are not Poissonian for any $\alpha$.  
\end{conj}
In \cite{not17} A.\ Walker proved for $(a_n) = (p_n)$ the sequence of primes that for almost all $\alpha$ the pair correlations of $(\lbrace p_n \alpha \rbrace)_{n \in \NN}$ are not Poissonian. Our conjecture would imply that there is no $\alpha$ such that $(\lbrace p_n \alpha \rbrace)_{n \in \NN}$ is Poissonian. \\

To be able to prove our result we need an alternative classification of integer sequences $\left(a_{n}\right)_{n \in \NN}$ with $E \left(A_{N}\right) = \Omega \left(N^{3}\right)$:\\

For $v \in \mathbb{Z}$ let $A_{N} (v)$ denote the cardinality of the set
$$
\left\{\left(x,y\right) \in \left\{1, \ldots, N\right\}^{2}, x \neq y : a_{x} - a_{y} = v \right\}.
$$
Then
\begin{equation} \label{equ_aa}
E \left(A_{N}\right) = \Omega \left(N^{3}\right)
\end{equation}
is equivalent to
\begin{equation} \label{equ_bb}
\sum_{v \in \mathbb{Z}} A^{2}_{N} (v) = \Omega \left(N^{3}\right),
\end{equation}
which implies that there is a $\kappa > 0$ and positive integers $N_{1} < N_{2} < N_{3} < \ldots$ such that
\begin{equation} \label{equ_c22}
\sum_{v \in \mathbb{Z}} A^{2}_{N_{i}} (v) \geq \kappa N^{3}_{i}, \qquad i = 1,2,\dots.~\\
\end{equation}

It will turn out that sequences $\left(a_{n}\right)_{n \in \NN}$ satisfying \eqref{equ_aa} have a strong linear substructure. From \eqref{equ_c22} we can deduce by the Balog--Szemeredi--Gowers-Theorem (see \cite{not5, not7}) that there exist constants $c, C > 0$ depending only on $\kappa$ such that for all $i=1,2,3,\ldots$ there is a subset $A_{0}^{(i)} \subset \left(a_{n}\right)_{1 \leq n \leq N_{i}}$ such that 
$$
\left|A_{0}^{(i)}\right| \geq c N_{i} \qquad \text{and} \qquad \left|A_{0}^{(i)} + A_{0}^{(i)}\right| \leq C  \left|A_{0}^{(i)}\right| \leq C N_{i}.
$$
The converse is also true: If for all $i$ for a set $A_{0}^{(i)}$ with $A_{0}^{(i)} \subset \left(a_{n}\right)_{1 \leq n \leq N_{i}}$ with $\left| A_{0}^{(i)} \right| \geq c N_{i}$ we have $\left|A_{0}^{(i)} + A_{0}^{(i)}\right| \leq C \left|A_{0}^{(i)}\right|$, then 
$$
\sum_{v \in \mathbb{Z}} A^{2}_{N_{i}} (v) \geq \frac{1}{C} \left|A_{0}^{(i)}\right|^{3} \geq \frac{c^3}{C} N_{i}^{3}
$$
and consequently $\sum_{v \in \mathbb{Z}} A^{2}_{N} (v)= \Omega \left(N^{3}\right)$ (this an elementary fact, see for example Lemma~1~(iii) in \cite{not10}.)\\

Consider now a subset $A_{0}^{(i)}$ of $\left(a_{n}\right)_{1 \leq n \leq N_{i}}$ with 
$$
\left|A_{0}^{(i)}\right| \geq c N_{i} \qquad \text{and} \qquad \left|A_{0}^{(i)} + A_{0}^{(i)} \right| \leq C \left|A_{0}^{(i)}\right|.
$$
By the theorem of Freiman (see \cite{not6}) there exist constants $d$ and $K$ depending only on $c$ and $C$, i.e., depending only on $\kappa$ in our setting, such that there exists a \emph{$d$-dimensional arithmetic progression} $P_{i}$ of size at most $K N_{i}$ such that $A_{0}^{(i)} \subset P_{i}$. This means that $P_i$ is a set of the form 
\begin{equation}\label{equ_c}
P_{i} := \left\{ \left. h_{i} + \sum^{d}_{j=1} r_{j} k_{j}^{(i)} \right|  0 \leq r_{j} < s_{j}^{(i)} \right\},  
\end{equation}
with $h_{i}, k_{1}^{(i)}, \ldots, k^{(i)}_{d}, s_{1}^{(i)}, \ldots, s_{d}^{(i)} \in \mathbb{Z}$ and such that $s_{1}^{(i)} s_{2}^{(i)} \ldots s_{d}^{(i)} \leq K N_{i}$.\\

In the other direction again it is easy to see that for any such set $A_{0}^{(i)}$
$$
\left|A_{0}^{(i)} + A_{0}^{(i)}\right| \leq 2^{d} K N_{i}.
$$
Based on these observations we make the following definition:
\begin{defi} \label{def_a}
Let $\left(a_{n}\right)_{n \in \NN}$ be a strictly increasing sequence of positive integers. We call this sequence {\em quasi-arithmetic of degree} $\mathbf{d}$, where $d$ is a positive integer, if there exist constants $C,K > 0$ and a strictly increasing sequence $\left(N_{i}\right)_{i \geq 1}$ of positive integers such that for all $i \geq 1$ there is a subset $A^{(i)} \subset \left(a_{n}\right)_{1 \leq n \leq N_{i}}$ with $\left|A^{(i)}\right| \geq C N_{i}$ such that $A^{(i)}$ is contained in a $d$-dimensional arithmetic progression $P^{(i)}$ of size at most $K N_{i}$.
\end{defi}
The above considerations show:

\begin{prop}
For a strictly increasing sequence $\left(a_{n}\right)_{n \in \NN}$ of positive integers we have $E\left(A_{N}\right)= \Omega \left(N^{3}\right)$ if and only if $\left(a_{n}\right)_{n \in \NN}$ is quasi-arithmetic of some degree $d$.
\end{prop}
Hence, our Theorem 1 stated above is equivalent to:
\begin{prop} \label{th_quasi}
If $\left(a_{n}\right)_{n \in \NN}$ is quasi-arithmetic of degree $d$, then there is \textbf{no} $\alpha$ such that the pair correlations of $\left(\left\{a_{n} \alpha\right\}\right)_{n \in \NN}$ are Poissonian.
\end{prop}
The result was already proven by the first author for $d=1$ in a previous work, see \cite{not18}. This case can also be recovered by Theorem 1 in \cite{not1}. 
\section{Proof of Theorem 1}
As noted above it is sufficient to prove Proposition 2. Let now $(a_n)_{n \in \NN}$ be quasi-arithmetic of degree $d$. That means (see Definition 2): There exists a strictly increasing subsequence $(N_i)_{i \in \NN}$ of the positive integers and $C,K>0$ with the following property: For all $i \geq 1$ there is a subset $b_1 < b_2 < \ldots < b_{M_i}$ of $(a_n)_{n=1, \ldots, N_i}$ with $M_i \geq CN_i$, such that $(b_j)_{j=1, \ldots, M_i}$ is a subset of 
\begin{equation*}
P_i:= \left\lbrace \left. h_i + \sum_{j=1}^d r_j k_j^{(i)} \right| \ 0 \leq r_j < s_j^{(i)} \right\rbrace
\end{equation*}
with certain $h_i, k_1^{(i)}, \ldots, k_d^{(i)} \in \ZZ$, $s_1^{(i)}, \ldots, s_d^{(i)} \in \NN$ and $s_1^{(i)}s_2^{(i)} \ldots s_d^{(i)} \leq KN_i$. Fix now any $i$, and for simplicity we omit the index $i$ in the above notations, i.e., we put $M:=M_i, h:=h_i$ and so on. In the sequel, we will put $K=1$ and $h=0$. The general case is treated similarly. Further, for $k=1, \ldots, M$, we set 
\begin{equation*}
b_k = r_1^{(k)} k_1 + \ldots + r_d^{(k)}k_d
\end{equation*}
and we identify $b_k$ with the vector
\begin{equation*}
(r_1^{(k)}, \ldots, r_d^{(k)})=:\boldsymbol{r}_k.
\end{equation*}
We have $0 \leq r_j^{(k)} < s_j$ for all $k=1, \ldots, M$ and all $j=1, \ldots, d$. Consider the differences $\boldsymbol{r}_k  - \boldsymbol{r}_l$ for $k,l=1, \ldots, M$. This yields $M^2 \geq C^2N^2$ vectors (counted with multiplicity) 
\begin{equation*}
\boldsymbol{u}:= \left( \begin{array}{c} u_1 \\ \vdots \\ u_d \end{array} \right),
\end{equation*}
with $-(s_j-1) \leq u_j \leq (s_j-1)$ for $j=1, \ldots, d$. There exist at most $2^d s_1 \ldots s_d \leq 2^d N$ \textit{different} such vectors $\boldsymbol{u}$. For each such given vector $\boldsymbol{u}$ there exist at most $M \leq N$ pairs $\boldsymbol{r}_k, \boldsymbol{r}_l$ such that $\boldsymbol{r}_k - \boldsymbol{r}_l = \boldsymbol{u}$. Let $\gamma:= \frac{C^2}{1+2^d}$, then there exist at least $\gamma N$ different vectors $\boldsymbol{u}$ such that there exist at least $\gamma N$ pairs $\boldsymbol{r}_k, \boldsymbol{r}_l$ with $\boldsymbol{r}_k - \boldsymbol{r}_l = \boldsymbol{u}$. Otherwise we had: 
\begin{align*}
C^2N^2  \leq M^2 &\leq \gamma N M +(2^dN -\gamma N)\gamma N \\
&\leq \gamma N^2 + (2^d-\gamma)\gamma N^2, 
\end{align*}
hence 
\begin{equation*}
C^2 \leq \gamma  + (2^d-\gamma)\gamma < \gamma (1+2^d),
\end{equation*}
i.e., $\gamma > \frac{C^2}{1+2^d}$, a contradiction. \\ \\
In the sequel, we will refer to this observation as \textit{Property 1}. Take now $\gamma N$ such $d$-tuples $\boldsymbol{u}$ having Property 1 and consider the corresponding $\gamma N$ values 
\begin{equation}\label{eq:eq1}
\lbrace (u_1 k_1 + \ldots + u_d k_d) \alpha \rbrace, \qquad \text{ in } [0,1). 
\end{equation}  
Let $L:= \frac{2}{\gamma}$, then there is a $\beta \in [0,1)$, such that the interval $\left[\beta, \beta + \frac{L}{\gamma N} \right)$ contains at least $L$ elements of the form (\ref{eq:eq1}), say the elements 
\begin{equation*}
\lbrace (u_1^{x} k_1 + \ldots + u_d^{x} k_d) \alpha \rbrace, \qquad \text{ for } x=1, \ldots, L. 
\end{equation*}
We call this fact \textit{Property 2}. \\

For every choice of $x$, we now consider $\gamma N$ pairs of $d$-tuples, say 
\begin{equation*}
\boldsymbol{r}_{i,x} := \left( \begin{array}{c} r_{1}^{(i, x)} \\ \vdots \\ r_{d}^{(i, x)} \end{array} \right), \text { and } \tilde{\boldsymbol{r}}_{i,x} :=\left( \begin{array}{c} \tilde{r}_{1}^{(i, x)} \\ \vdots \\ \tilde{r}_{d}^{(i, x)} \end{array} \right),
\end{equation*}
for $i=1, \ldots, \gamma N$, such that 
\begin{equation*}
\boldsymbol{r}_{i,x} - \tilde{\boldsymbol{r}}_{i,x} =  \left( \begin{array}{c} u_1^{(x)} \\ \vdots \\ u_d^{(x)} \end{array} \right).
\end{equation*}
We will show that there exist $x, y \in \lbrace 1, \ldots, L \rbrace$ with $x \neq y$ such that 
\begin{equation}\label{eq:eq2}
\# \left( \lbrace \boldsymbol{r}_{i,x} | i=1, \ldots, \gamma N \rbrace \cap \lbrace \boldsymbol{r}_{i,y} | i = 1, \ldots, \gamma N \rbrace \right) \geq \frac{N}{L^2}. 
\end{equation} 
Assume this were not the case and define 
\begin{equation*}
\mathcal{M}_x:= \lbrace \boldsymbol{r}_{i,x} | i=1, \ldots, \gamma N \rbrace.
\end{equation*}
Then, we had 
\begin{align*}
N \geq s_1 \ldots s_d &\geq \Big| \bigcup_{x=1}^{L} \mathcal{M}_x \Big| \\
& \geq \sum_{x=1}^L | \mathcal{M}_x |  - \sum_{\substack{x,y=1 \\ x \neq y}}^L | \mathcal{M}_x \cap \mathcal{M}_y | \\
& > L \gamma N - L^2 \frac{N}{L^2} = N, 
\end{align*}
which is a contradiction. \\

Let now $x$ and $y$ satisfying (\ref{eq:eq2}) be given. Let 
\begin{align*}
\boldsymbol{r}_i&, \quad i=1, \ldots, \frac{N}{L^2}, \\
\tilde{\boldsymbol{r}}_{i,x}&, \quad i=1, \ldots, \frac{N}{L^2}, \\
\tilde{\boldsymbol{r}}_{i,y}&, \quad i=1, \ldots, \frac{N}{L^2} 
\end{align*} 
be such that 
\begin{equation*}
\boldsymbol{r}_i - \tilde{\boldsymbol{r}}_{i,x} = \left( \begin{array}{c} u_1^{(x)} \\ \vdots \\ u_d^{(x)} \end{array} \right), \text{ and } 
\boldsymbol{r}_i - \tilde{\boldsymbol{r}}_{i,y} = \left( \begin{array}{c} u_1^{(y)} \\ \vdots \\ u_d^{(y)} \end{array} \right).
\end{equation*}
Then,
\begin{equation*}
\tilde{\boldsymbol{r}}_{i,y}-\tilde{\boldsymbol{r}}_{i,x} = \left( \begin{array}{c} u_1^{(x)} - u_1^{(y)} \\ \vdots \\ u_d^{(x)} -u_d^{(y)} \end{array} \right) =: \left( \begin{array}{c} z_1 \\ \vdots \\ z_d \end{array} \right), 
\end{equation*}
for $i=1, \ldots, \frac{N}{L^2}$. Due to Property 2, we have
\begin{equation*}
\frac{L}{\gamma N} \geq | \lbrace (z_1k_1 + \ldots + z_d k_d) \alpha \rbrace |.
\end{equation*}
To sum up, we have shown that for all $N_i$ there exist at least 
\begin{equation*}
\frac{N_i}{L^2} = \frac{4(1+2^d)^2}{C^4}N_i =: \tau N_i,
\end{equation*}
pairs $(k,l)$ with $1 \leq k \neq l \leq N_i$, such that all expressions $\| \lbrace a_k \alpha \rbrace - \lbrace a_l \alpha \rbrace \|$ have the same value and satisfy 
\begin{equation*}
\| \lbrace a_k \alpha \rbrace - \lbrace a_l \alpha \rbrace \| \leq \frac{L}{\gamma N_i} = \frac{2(1+2^d)^2}{C^4} \frac{1}{N_i} =: \psi \frac{1}{N_i}.
\end{equation*}
Note, that $\tau$ and $\psi$ only depend on $d$ and $C$ (and on $K$ if $K \neq 1$). For every $i$ choose now $\psi_i$ minimal such that there exist at least $\tau N_i$ pairs $(k,l)$ with $1 \leq k \neq l \leq N_i$, such that 
\begin{equation*}
\| \lbrace a_k \alpha \rbrace - \lbrace a_l \alpha \rbrace \| = \psi_i \frac{1}{N_i}. 
\end{equation*}  
Of course, $\psi_i \leq \psi$ for all $i$. Let now $\rho:= \frac{\tau}{3}$ and assume that $\psi_i < \rho$ for infinitely many $i$. Therefore, we have for these $i$
\begin{align*}
\frac{1}{N_i} &\# \left\lbrace 1 \leq k \neq l \leq N_i | \ \| \lbrace a_k \alpha \rbrace - \lbrace a_l \alpha \rbrace \| \leq \rho \frac{1}{N_i} \right\rbrace \\
& \geq \tau = 3 \rho, 
\end{align*}
which is a contradiction and consequently the pair correlations are not Poissonian. \\

Assume now that $\psi_i \geq \rho$ for infinitely many $i$. Consequently, there exists an $s_1 \in [\rho, \psi)$ such that 
\begin{equation*}
\psi_i \in \left[s_1, s_1 + \frac{\tau}{3} \right)
\end{equation*}
for infinitely many $i$. In the following, we only consider these $i$ and we will set $s_2:=  s_1 + \frac{2\tau}{3}$. Then, we have
\begin{align*}
&\frac{1}{N_i} \# \left\lbrace 1 \leq k \neq l \leq N_i | \ \| \lbrace a_k \alpha \rbrace - \lbrace a_l \alpha \rbrace \| \leq s_2 \frac{1}{N_i} \right\rbrace \\
&- \frac{1}{N_i} \# \left\lbrace 1 \leq k \neq l \leq N_i | \ \| \lbrace a_k \alpha \rbrace - \lbrace a_l \alpha \rbrace \| \leq s_1 \frac{1}{N_i} \right\rbrace \\
& \geq \tau.
\end{align*}
If $(\lbrace a_n \alpha \rbrace)_{n \in \NN}$ were Poissonian, then the above difference should converge, as $i \to \infty$, to $2(s_2 -s_1)= \frac{2\tau}{3} < \tau$, which is a contradiction. 
\noindent \\ \\
\textbf{Acknowledgements:} The authors thank an anonymous referee for pointing out an inaccuracy in the first version of the paper. 
  
\textbf{Author's Addresses:} \\ 
Gerhard Larcher and Wolfgang Stockinger, Institut f\"ur Finanzmathematik und Angewandte Zahlentheorie, Johannes Kepler Universit\"at Linz, Altenbergerstra{\ss}e 69, A-4040 Linz, Austria. \\ \\ 
Email: gerhard.larcher(at)jku.at, wolfgang.stockinger(at)jku.at       

\begin{thebibliography}{99}
\addcontentsline{toc}{section}{References}
\bibitem{not2} I.\ Aichinger, C.\ Aistleitner, G.\ Larcher, \textit{On Quasi-Energy-Spectra, Pair Correlations of Sequences and Additive Combinatorics},
Celebration of the 80th birthday of Ian Sloan (J.\ Dick, F.\ Y.\ Kuo, H.\ Wo\'zniakowski, eds.), Springer-Verlag, to appear, 2018. 
\bibitem{not3} C.\ Aistleitner, G.\ Larcher, M.\ Lewko, \textit{Additive energy and the Hausdorff dimension of the exceptional set in metric pair correlation problems. With an appendix by Jean Bourgain}, Israel J.\ Math., 222 (2017) No.\ 1, $463$-$485$. 
\bibitem{not4} C.\ Aistleitner, T.\ Lachmann, F.\ Pausinger, \textit{Pair correlations and equidistribution}, Journal of Number Theory, to appear, available at \url{https://arxiv.org/abs/1612.05495}.
\bibitem{not5} A.\ Balog, E.\ Szemer\'edi, \textit{A statistical theorem of set addition}, Combinatorica, 14:263--268, 1994.
\bibitem{not6} G.\ R.\ Freiman, \textit{Foundations of a Structural Theory of Set Addition}, Translations of Mathematical Monographs 37, Amer.\ Math.\ Soc., Providence, USA, 1973.
\bibitem{not7} W.\ T.\ Gowers, \textit{A new proof of Szemer\'edi's theorem for arithmetic progressions of length four}, Geometric and Functional Analysis, 8:529--551, 1998.
\bibitem{not8} S.\ Grepstad, G.\ Larcher, \textit{On pair correlation and discrepancy}, Arch.\ Math., 109:143--149, 2017.
\bibitem{not9} T.\ Lachmann, N.\ Technau, \textit{On Exceptional Sets in the Metric Poissonian Pair Correlations problem}, Submitted, 2017. arXiv:1708.08599.
\bibitem{not18} G.\ Larcher, \textit{Remark on a result of Bourgain on Poissonian pair correlation}, 
arXiv:1711.08663, 2017.
\bibitem{not1} G.\ Larcher, W.\ Stockinger, \textit{Some negative results related to Poissonian pair correlation problems}, arXiv:1803.05236, 2018.
\bibitem{not10} V.\ F.\ Lev, \textit{The (Gowers--)Balog--Szemer\'edi Theorem: An Exposition}, \url{http://people.math.gatech.edu/~ecroot/8803/baloszem.pdf}. 
\bibitem{not11} J.\ Marklof, \textit{The Berry--Tabor conjecture}, European Congress of Mathematics, Vol. II (Barcelona, 2000), 421--427, Progr.\ Math., 202, Birkh\"auser, Basel, 2001.
\bibitem{not12} Z.\ Rudnick, P.\ Sarnak, \textit{The pair correlation function of fractional parts of polynomials}, Comm.\ Math.\ Phys., 194(1):61--70, 1998.
\bibitem{not13} Z.\ Rudnick, P.\ Sarnak, A.\ Zaharescu, \textit{The distribution of spacings between the fractional parts of $n^{2} \alpha$}, Invent.\ Math., 145(1):37--57, 2001.
\bibitem{not14} Z.\ Rudnick, A.\ Zaharescu, \textit{A metric result on the pair correlation of fractional parts of sequences}, Acta Arith., 89(3):283--293, 1999.
\bibitem{not15} S.\ Steinerberger, \textit{Localized Quantitative Criteria for Equidistribution}, Acta Arith., 180:183--199, 2017.
\bibitem{not16} T.\ Tao, V.\ Vu, \textit{Additive combinatorics}, Volume 105 of Cambridge Studies in Advanced Mathematics, Cambridge University Press, Cambridge, 2006.
\bibitem{not17} A.\ Walker, \textit{The primes are not Poissonian}, arXiv:1702.07365, 2017. 
\end{thebibliography}
\end{document}